\documentclass[12pt]{amsproc}
\usepackage[cp1251]{inputenc}
\usepackage[russian]{babel}
\usepackage{amssymb}
\numberwithin{equation}{section}
\newtheorem{tm}{Теорема}[section]
\newtheorem*{tmP}{Теорема Понтрягина}
\newtheorem*{tmL}{Теорема Лангера}
\newtheorem*{tmKLA}{Теорема Крейна-Лангера-Азизова}
\newtheorem*{tmOG}{Теорема о генераторе \(C_0\)-полугруппы}
\newtheorem*{tmOGG}{Теорема о генераторе голоморфной полугруппы}
\newtheorem*{tmGo}{Теорема Гомилко}
\newtheorem*{tmGe}{Теорема Gearhart'а}
\newtheorem{lem}[tm]{Лемма}
\newtheorem{rem}[tm]{Замечание}

\newcommand{\rank}{\operatorname{rank}}
\renewcommand{\Re}{\operatorname{Re}}
\renewcommand{\Im}{\operatorname{Im}}
\newcommand{\Lin}{\operatorname{Lin}}
\renewcommand{\ge}{\geqslant}
\renewcommand{\le}{\leqslant}
\newcommand{\la}{\lambda}

\newcommand{\ve}{\varepsilon}
\newcommand{\al}{\alpha}
\newcommand{\ov}{\overline}
\newcommand{\p}{^{-1}}

\title[]{Диссипативные операторы в
пространстве Крейна. Инвариантные подпространства и свойства
сужений}
\author[]{А.~А.~Шкаликов}
\thanks{Работа поддержана программой \flqq Ведущие научные
школы\frqq, грант \No~НШ-5247.2006.1, Российским фондом
фундаментальных исследований и Немецким научно-исследовательским
обществом~DFG}

\begin{document}
\begin{abstract}
В работе доказана теорема о существовании максимального
неотрицательного инвариантного подпространства у  диссипативного
оператора в пространстве Крейна, который допускает матричное
представление относительно канонической декомпозиции пространства,
причем правый верхний оператор в этой декомпозиции компактен
относительно правого нижнего. При дополнительном предположении
ограниченности левых верхнего и нижнего операторов (так
называемого условии Лангера) этот результат был доказан (в порядке
возрастания общности) Понтрягиным, Крейном, Лангером и Азизовым.
Условие Лангера заменено в работе существенно более слабым. При
его выполнении доказано, что максимальный диссипативный оператор в
пространстве Крейна обладает максимальным, неотрицательным,
инвариантным подпространством, таким, что спектр сужения оператора
на это подпространство лежит в левой полуплоскости. Найдены
достаточные условия для того, чтобы сужение оператора на это
подпространство было генератором голоморфной или $C_0$-полугруппы.
\end{abstract}
\begin{flushleft}
УДК 517.9+517.43
\end{flushleft}
\maketitle

\section{Введение}

Напомним основные определения теории пространств с индефинитной
метрикой. Пусть \(\mathcal H\) "--- гильбертово пространство, в
котором наряду с обычным скалярным произведением \((x,y)\) задано
индефинитное скалярное произведение \([x,y]=(Jx,y)\), где $J$ "---
самосопряженный  оператор, такой, что $J^2 =I.$ Здесь и далее
через $I$ обозначается тождественный оператор, а пространство
\(\mathcal H\) считается сепарабельным. Оператор $J$ называют
канонической или внутренней симметрией в \(\mathcal H\). При
изучении гамильтоновых систем роль $J$   играет симплектическая
матрица. Поскольку оператор $J$ самосопряженный и унитарный
одновременно, то его спектр состоит из двух точек $\{\pm 1\}$.
Поэтому  справедливо представление \(J=P_+-P-\), где \(P_+\) и
\(P_-\) "--- взаимно ортогональные ортопроекторы, причем
\(P_++P_-=I\). Пространство
$(\, \mathcal H, [\, \cdot ,\cdot\, ]\, )$ называют пространством
Понтрягина и обозначают
\(\Pi_{\varkappa}\), если одно из чисел \(\rank P_+\) или \(\rank
P_-\) конечно и равно
\(\varkappa\), и пространством Крейна, если оба этих числа
бесконечны. Подпространство \(\mathcal L\) в $(\, \mathcal H, [\,
\cdot ,\cdot\, ]\, )$ называют
неотрицательным (положительным), если для всех
\(x\in\mathcal L\) выполнено неравенство $ [x,x]\geqslant 0 \
\, (>0)$. Если неотрицательное (положительное) подпространство
\(\mathcal L\) не допускает
неотрицательных (положительных)нетривиальных расширений, то оно
называется максимальным. Аналогично определяются максимальные
неположительные и отрицательные подпространства в $(\,
\mathcal H, [\, \cdot ,\cdot\, ]\, )$.

Пусть \(A\) "--- линейный оператор в \(\mathcal H\) с областью
определения \(\mathcal D(A)\). Через \(\sigma(A)\) и \(\rho(A)\)
обозначаем его спектр и резольвенту, соответственно. Оператор \(A\)
называется {\it диссипативным в} \(\mathcal H\), если
\(\Re(Ax,x)\leqslant 0\) для всех \(x\in\mathcal D(A)\).
Диссипативный оператор называется {\it максимально диссипативным}
(или {\it m-диссипативным}), если он не допускает нетривиальных
диссипативных расширений. Известно~\cite[Гл.~V, \S\,3.10]{1}, что
последнее условие для диссипативного оператора \(A\) эквивалентно
условию
\(\rho(A)\supset\mathbb C^-\), где \(\mathbb C^-\) "--- открытая
левая комплексная  полуплоскость. Оператор \(A\) называется
\(J\)-{\it диссипативным} (\(J\)-{\it максимально диссипативным}),
{\it или диссипативным (максимально диссипативным) в пространстве}
$(\, \mathcal H, [\cdot\, ,\cdot\, ]\, )$, если \(JA\) "--- диссипативный
(максимально диссипативный) в \(\mathcal H\). Аналогично, \(A\)
называется {\it симметрическим (самосопряжённым) в пространстве}
$(\, \mathcal H, [\cdot\, ,\cdot\, ]\, )$, если \(JA\) симметричен
(самосопряжён) в пространстве \(\mathcal H\).

Представим пространство \(\mathcal H\) в виде ортогональной суммы
\(\mathcal H=\mathcal H^+\oplus\mathcal H^-\), где \(\mathcal
H^{\pm}=P_{\pm}\mathcal H\) "--- образы проекторов \(P_{\pm}\). Пусть
\(A\) "--- линейный оператор в \(\mathcal H\) с областью определения
\(\mathcal D(A)\). Обозначим  \(\mathcal D^{\pm}=
\mathcal D(A)\cap \mathcal H^{\pm}\). Далее мы будем иметь
дело только с такими операторами, для которых сумма \(\mathcal D^+
\oplus\mathcal D^-\) является {\it ядром оператора} $A$. Это означает
следующее: если $A_0$ "--- сужение оператора $A$  на \(\mathcal
D^+ \oplus\mathcal D^-\), то $A \subset \bar{A_0}$, где
$\bar{A_0}$ "--- замыкание $A$.

Рассматривая, если нужно, сужение \(A\), будем далее считать, что
\(\mathcal D(A)=\mathcal D^+\oplus\mathcal D^-\). В этом случае
\(A\) допускает представление в виде операторной матрицы по
отношению к канонической декомпозиции \(\mathcal H=\mathcal H^+
\oplus\mathcal H^-\):
\begin{equation}\label{eq:1}
    A=\begin{pmatrix}P_+AP_+&P_+AP_-\\ P_-AP_+&
    P_-AP_-\end{pmatrix}:=\begin{pmatrix}A_{11}&A_{12}\\
    A_{21}&A_{22}\end{pmatrix}.
\end{equation}
В таком представлении векторы \(x=x_++x_-\in\mathcal H\), где
\(x_{\pm}\in\mathcal H^{\pm}\), отождествляются со столбцами
\(x=\begin{pmatrix}x_+\\ x_-\end{pmatrix}\), и действие \(A\)
задаётся равенством
\[
    Ax=A\begin{pmatrix}x_+\\ x_-\end{pmatrix}=
    \begin{pmatrix} A_{11}x_++A_{12}x_-\\
    A_{21}x_++A_{22}x_-\end{pmatrix},\qquad
    x_+\in\mathcal D^+,\; x_-\in\mathcal D^-.
\]

В 1941 году С.Л.Соболев обратил внимание Л.С.Понтрягина на
следующий замечательный факт: {\it самосопряженный оператор в
прстранстве $\Pi_1$ всегда имеет по меньшей мере один собственный
вектор} (доказательство Соболева долгое время оставалось
неопубликованным; этот результат  Соболев применял для
исследования устойчивости одной задачи механики \cite {1a}). Это
наблюдение Соболева послужило толчком для исследования Понтрягина
~\cite{2}, где он заложил начала геометрии пространств с
индефинитной метрикой и получил следующий фундаментальный
результат.

\begin{tmP}
Пусть \(A\) "--- самосопряжённый оператор в пространстве
\(\Pi_{\kappa} \), причём \(\rank P_+=\varkappa<\infty\). Тогда
существует \(A\)-инвариантное максимальное неотрицательное
подпространство \(\mathcal L\) (\(\dim\mathcal L=\varkappa\))
такое, что спектр сужения \(A|_{\mathcal L}\) лежит в замкнутой
верхней полуплоскости.
\end{tmP}

После  работы~\cite{2}, проблема существования инвариантных
максимальных полудефинитных подпространств находилась в центре
внимания теории операторов в пространствах Понтрягина и Крейна.
Крейн~\cite{2a} получил аналог теоремы Понтрягина для
нерастягивающих операторов в \(\Pi_{\varkappa}\), развивая другой
подход к проблеме, нежели в~\cite{2}. Важное обобщение теоремы
Понтрягина было получено Лангером~\cite{3,4} и Крейном~\cite{5}.
Сформулируем результат~\cite{4}.

\begin{tmL}
Пусть \(A\) "--- самосопряжённый оператор в пространстве Крейна
\((\mathcal H,[\cdot\, ,\cdot\, ])\), причём \(\mathcal D(A)\supset\mathcal H^+\)
(это условие эквивалентно возможности представления~\eqref{eq:1},
в котором операторы \(A_{11}\) и \(A_{12}\) ограничены).
Если оператор \(A_{12}=P_+AP_-\) компактен, то существует
\(A\)-инвариантное максимальное неотрицательное подпространство
\(\mathcal L\) такое, что спектр сужения \(A|_{\mathcal L}\)
лежит в замкнутой верхней полуплоскости.
\end{tmL}

Впоследствии теоремы о существовании \(A\)-инвариантных
подпространств были получены для других классов операторов.
Крейн ввёл и исследовал класс дефинитных операторов, а позже
для более широкого класса так называемых дефинизируемых операторов
Лангер~\cite{6,7} доказал теорему о существовании максимальных
дефинитных инвариантных подпространств, а также аналог
спектральной теоремы для самосопряжённых операторов.
 Крейн и Лангер~\cite{8}, и независимо Азизов~\cite{9}, начали исследование
 диссипативных операторов в пространствах Понтрягина и Крейна.
В частности, в  этих работах они доказали следующий результат.

\begin{tmKLA}
Пусть $-iA$ "--- максимально диссипативный оператор в пространстве
Понтрягина \(\Pi_{\kappa} \). Тогда справедливо утверждение
теоремы Понтрягина.
\end{tmKLA}

Далее Азизов и Хорошавин~\cite{10} доказали аналог теоремы Лангера
для некоторого класса нерастягивающих операторов в пространстве
Крейна, а с помощью этого результата Азизов~\cite[Гл.~2]{11}
доказал, что теорема Лангера~\cite{4} остаётся справедливой для
максимально диссипативных операторов в пространстве Крейна (при
этом в формулировке  верхнюю комплексную полуплоскость нужно
заменить на левую). Более того, Азизов показал, что условие
компактности оператора $A_{12}$ можно заменить на условие его
компактности по отношению к оператору $A_{22}$, т.е. компактности
$A_{12}(A_{22} -\mu)^{-1}$ при некотором $\mu\in\rho(A_{22})$ (для
самосопряженных операторов этот факт ранее был отмечен Крейном
\cite{5}.) Прямое, более короткое доказательство теоремы Азизова
было предложено автором~\cite{12}.

Задачу описания вещественного спектра или чисто мнимого спектра
для сужений самосопряжённых или диссипативных операторов в пространстве
Крейна на максимальное неотрицательное инвариантное подпространство
изучали Костюченко и Оразов~\cite{17}, Гомилко~\cite{16}
и автор~\cite{19}.

{\it Условие Лангера} \(\mathcal D(A)\supset\mathcal H^+\) (или
условие ограниченности в матричном представлении~\eqref{eq:1}
операторов \(A_{11}\) и \(A_{21}\)) является весьма
ограничительным. В частности, оно не выполняется в некоторых
конкретных задачах, которые были модельными в
работах~\cite{14,15}). Это обстоятельство послужило мотивом для
получения  обобщений сформулированых результатов, когда условие
Лангера заменяется другим более слабым. Такой результат был
недавно получен автором \cite{18a}. Условие Лангера в \cite{18a}
заменено следующим условием: найдётся число \(\mu\) в левой полуплоскости
такое, что операторы
\[
    (A_{22}-\mu)^{-1}A_{21},\qquad A_{21}(A_{22}-\mu)^{-1},\qquad
    S(\mu)=A_{11}-A_{12}(A_{22}-\mu)^{-1}A_{21}
\]
ограничены. \(S(\mu)\) называют {\em передаточной функцией} (а
также трансфер-функцией или дополнением Фробениуса--Шура). Язык,
использующий передаточные функции, на первый взгляд может
показаться  неуклюжим, но он естествен и соответствующие условия в
конкретных задачах нетрудно проверять. При этом условие
ограниченности первых двух операторов, как правило, выполняется
автоматически. Но условие ограниченности передаточной функции
является достаточно  "<жестким">. Ниже мы увидим, что оно
выполняется в том и только в том случае, когда сужение оператора
$A$  на соответствующее инвариантное подпространство является
ограниченным оператором. Конечно, имея в виду дальнейшие
приложения (которые в этой работе мы не рассматриваем),
представляется привлекательной задача от этого условия избавиться.
{\it Это является основной целью настоящей работы}. Но научившись
находить инвариантные подпространства, такие, что сужения
$A|_{\mathcal L}$ неограничены, мы приходим к {\em новой важной
задаче}: являются ли они генераторами какого-либо типа полугрупп?
{\it Найти достаточные условия для положительного ответа на этот
вопрос --- вторая цель работы.}

\section{Существование инвариантных подпространств}

В цитированных выше работах предполагалось, что подпространство
$\mathcal L$ является $A$-инвариантным, если $\mathcal L \subset\mathcal
D(A)$ и $A(\mathcal L) \subset \mathcal L$. Но этом случае сужение
$A|_\mathcal L$ ограниченный оператор. Поэтому далее считаем, что
подпространство $\mathcal L$ является $A-$ инвариантным, если
$\mathcal D(A)\cap\mathcal L$ плотно в $\mathcal L$
 и $Ax \in \mathcal L$ для всех $x\in \mathcal D(A) \cap\mathcal
 L$.

 Основная теорема об
инвариантных подпространствах будет доказана для максимально
диссипативных операторов в пространстве Крейна.
Известно~\cite[Гл.~2, Теорема~2.9]{11}, что при выполнении условия
Лангера \(\mathcal D(A)\supset\mathcal H^+\) оператор \(A\)
максимально диссипативен в $(\,\mathcal H, [\cdot\ ,\cdot ]\, )$
 тогда и только тогда, когда \(-A_{22}\)
максимально диссипативен в \(\mathcal H^-\), т.е. левая
полуплоскость лежит в резольвентном множестве оператора
\(A_{22}\). В общем случае можно утверждать, что максимальная
диссипативность оператора $A$ влечет только диссипативность
операторов $A_{11}$ и $-A_{22}$. Но мы будем работать с
передаточной функцией, поэтому будем требовать условие
существования резольвенты при некотором (а тогда при всех)
$\mu\in\mathbb C^-$. Сущность налагаемых ниже условий состоит в том,
что оператор \(A_{22}\) является доминирующим по отношению к
сплетающим операторам \(A_{21}\) и \(A_{12}\).
 Сформулируем основной результат.

\begin{tm}\label{tm:1}
Пусть \(A\) "--- диссипативный оператор в пространстве Крейна $(\,
\mathcal H, [\cdot\, ,\cdot\, ]\, )$ с областью \(\mathcal D(A)=\mathcal D^+
\oplus\mathcal D^-\), плотной в \(\mathcal H=\mathcal H^+
\oplus\mathcal H^-\). Обозначим через $\overline A$ его замыкание.
Пусть~\eqref{eq:1} "--- матричное представление \(A\) в \(\mathcal
H^+\oplus\mathcal H^-\), и пусть выполнены следующие условия:
\begin{itemize}
\item[(i)] оператор \(-A_{22}\) является максимально диссипативным
в пространстве \(\mathcal H^-\) (т.е. существует резольвента
\((A_{22}- \mu)^{-1}\) при всех \(\mu\in\mathbb C^-\));
\item[(ii)] оператор \(F(\mu)=(A_{22}-\mu)^{-1}A_{21}\) допускает
ограниченное замыкание при некотором (а тогда при всех) \(\mu\in
\mathbb C^-\);
\item[(iii)] оператор \(G(\mu)=A_{12}(A_{22}-\mu)^{-1}\) компактен
при некотором (а тогда при всех) \(\mu\in\mathbb C^-\).
\end{itemize}

Тогда существует максимальное неотрицательное подпространство,
инвариантное относительно оператора  \(\overline{A}\). Если
оператор \(\overline A\) является \(m\)-диссипативным в $(\, \mathcal H,
[\cdot\, ,\cdot\,  ]\, )$ (и только в этом случае), то найдется максимальное
неотрицательное \(\overline{A}\)-инвариантное подпространство \(\mathcal L^+\)
такое, что спектр сужения \(A_+=\overline{A}|_{\mathcal L^+}\) содержится
в замкнутой левой полуплоскости.
\end{tm}

Доказательству теоремы предпошлём ряд лемм.  Первые два ---
известные предложения: первая лемма имеется в работе Понтрягина
\cite{2}, вторая доказана в \cite{14}. Основную роль в дальнейшем играют
Леммы~\ref{lem:4}, \ref{lem:5} и~\ref{lem:7}.

\begin{lem}\label{lem:1}
Подпространство \(\mathcal L\) является максимальным
неотрицательным  в \((\mathcal H,[\cdot\,,\cdot\,])\) в том и
только том случае, когда оно допускает представление
\begin{equation}\label{eq:2}
    \mathcal L=\left\{x=\begin{pmatrix}x_+\\ Kx_+\end{pmatrix},
    \qquad x_+\in\mathcal H^+\right\},
\end{equation}
где \(K:\mathcal H^+\to\mathcal H^-\) "--- линейный
нерастягивающий оператор (т.е.
\(\|K\|\leqslant 1\)).
\end{lem}

Оператор \(K\) в представлении~\eqref{eq:2} называется \emph{угловым
оператором} подпространства \(\mathcal L\).

\begin{lem}\label{lem:2}
Пусть \(A\) "--- оператор с плотной областью \(\mathcal D(A)=
\mathcal D^+\oplus\mathcal D^-\), резольвентное множество
\(\rho(A_{22})\) непусто, а операторы
\[
    G=A_{12}(A_{22}-\mu)^{-1},\qquad F=(A_{22}-\mu)^{-1}A_{21},
    \quad \text{и} \quad  S=A_{11}-A_{12}F,\qquad \mu\in\rho(A_{22}),
\]
ограничены при некотором \(\mu\in\rho(A_{22})\). Тогда оператор
\(A\) замыкаем в том и только том случае, когда \(S= A_{11} -A_{12}F(\mu)\) замыкаем
в \(\mathcal H^+\), и его замыкание определяется равенством
\begin{equation}\label{eq:4}
    \overline{A}=\mu+\begin{pmatrix}1&G\\0&1\end{pmatrix}
    \begin{pmatrix}\overline{S}-\mu&0\\0&A_{22}-\mu\end{pmatrix}
    \begin{pmatrix}1&0\\ \overline F&1\end{pmatrix}.
\end{equation}
Более точно, область определения и действие \(\overline{A}\)
определяются равенствами
\begin{gather*}
    \mathcal D(\overline{A})=\left\{\begin{pmatrix}x_+\\
    x_-\end{pmatrix}\in\mathcal H,\qquad x_+\in\mathcal H^+,
    \qquad \overline Fx_++x_-\in\mathcal D^-\subset\mathcal D(A_{22})
    \right\},\\
    \overline{A}\begin{pmatrix}x_+\\x_-\end{pmatrix}=
    \begin{pmatrix}\overline{S}x_++G(A_{22}-\mu)(\overline Fx_++x_-)\\
    (A_{22}-\mu)(\overline Fx_++x_-)+\mu x_-\end{pmatrix}.
\end{gather*}
\end{lem}

\begin{lem}\label{lem:3}
Пусть \(-A_{22}\) "--- максимально диссипативный, а
\(G(\mu)=A_{12}(A_{22}-\mu)^{-1}\) "--- компактный операторы при
некотором (а тогда при всех) \(\mu\in\mathbb C^-\). Тогда
\(\|G(\mu)\|\to 0\) при \(\mu\to \infty\) равномерно в любом секторе
$|\arg (\pi- \mu)|\leqslant \gamma <\pi/2$.
\end{lem}
\begin{proof}
Другая версия этого утверждения будет сформулирована в Лемме 3.6.
Там же будет приведено доказательство.
\end{proof}


\begin{lem}\label{lem:4}
При $\mu \in\rho(A_{22})$  справедливо равенство
\begin{equation}\label{eq:5}
    JA+\mu=J\begin{pmatrix}1&G\\0&1\end{pmatrix}
    \begin{pmatrix}S+\mu&0\\0&A_{22}-\mu\end{pmatrix}
    \begin{pmatrix}1&0\\F&1\end{pmatrix}  \qquad
    J:= \begin{pmatrix}1&0\\0&-1\end{pmatrix}.
\end{equation}
\end{lem}
\begin{proof}
В справедливости этого равенства легко убедиться непосредственной
проверкой с учетом представления~\eqref{eq:4}.
\end{proof}
\begin{lem}\label{lem:5}
При \(\mu\in\rho(A_{22})\) для всех
\(x_+\in\mathcal D_+\) справедливо равенство
\begin{equation}\label{eq:6}
    (Sx_+,x_+)=\left(JA\begin{pmatrix}x_+\\-Fx_+\end{pmatrix},
    \begin{pmatrix}x_+\\-Fx_+\end{pmatrix}\right)+\mu(Fx_+,Fx_+).
\end{equation}
В частности, при любом \(\mu\in\mathbb C^-\) оператор $S=S(\mu)$ c
с областью \(\mathcal D(S)=\mathcal D^+\) является диссипативным в
пространстве $\mathcal H^+$, если оператор $A$ диссипативен в
пространстве $(\, \mathcal H, [\, \cdot ,\cdot\, ]\, )$. Замыкание
оператора $S$ является $m$-диссипативным в пространстве $\mathcal
H^+$, если замыкание оператора $A$ является  $m$-диссипативным в
пространстве $(\, \mathcal H, [\,
\cdot ,\cdot\, ]\, )$, и наоборот.
\end{lem}
\begin{proof}
 В справедливости равенства \ref{eq:6} легко убедиться
непоспедственной проверкой с помощью равенства \ref{eq:5}.
Эквивалентность $m$-диссипативности операторов $S(\mu)$ в
пространстве $\mathcal H_+$ и $JA$  в пространстве $\mathcal H$
также следует из равенства \ref{eq:5}.
\end{proof}

\begin{lem}\label{lem:6}
Пусть подпространство  \(\mathcal L\) имеет вид
\[
    \mathcal L=\{x\;:\;x=\begin{pmatrix}x_+\\Kx_+\end{pmatrix},\;
    x_+\in \mathcal H^+\},
\]
где  \(K:\mathcal H^+\to \mathcal H^-\) --- ограниченный оператор.
Тогда
\(\mathcal L\) является  \(A\)-инвариантным в том и только том случае,
когда выполняется равенство
\begin{equation}\label{Ric1}
    (1-KG)(A_{22}-\mu)(F+K)=K(S-\mu).
\end{equation}

 Здесь предполагается, что областью определения оператора в
правой (а потому и в левой) части  равенства является линеал
\(\mathcal D(S)=\mathcal D^+\).
\end{lem}

\begin{proof}
Для  \(x_+\in\mathcal D_+\) имеем
\[
    (A-\mu)\begin{pmatrix}x_+\\Kx_+\end{pmatrix}=\begin{pmatrix}
    (S-\mu)x_++G(A_{22}-\mu)(F+K)x_+\\ (A_{22}-\mu)(F+K)x_+
    \end{pmatrix}.
\]
Если  \(\mathcal L\) является  \(A\)-инвариантным, то при
некотором \(y_+\in \mathcal H^+\) имеем

\begin{gather*}
    [(S-\mu)+G(A_{22}-\mu)(F+K)]x_+=y_+,\\
    (A_{22}-\mu)(F+K)x_+=Ky_+.
\end{gather*}
Подставляя выражение для $y_+$  во второе равенство, получаем
уравнение ~\eqref{Ric1}.

Обратно, если выполнено уравнение ~\eqref{Ric1}, то при некотором
$y_+$ выполнены два последних уравнения,  из которых получаем,
что подпространстсво \(\mathcal L\) является
\(A\)-инвариантным.
\end{proof}

\begin{rem}
Уравнение~\eqref{Ric1} является модифицированным уравнением
Риккати для неизвестного оператора $K$. Ранее, начиная с работы
Понтрягина~\cite{2}, использовалась другая версия
Леммы~\ref{lem:6}: \(\mathcal L\) является
\(A\)-инвариантным тогда и только тогда, когда

\begin{equation}\label{eq:R2}
    A_{21}+A_{22}K-KA_{11}-KA_{12}K=0.
\end{equation}
Это уравнение называют уравнением Риккати, ассоциированным с
оператор-матрицей $A$. Однако эта форма уравнения Риккати неудобна
в случае, когда все элементы оператор-матрицы $A$ могут быть
неограниченными операторами.
\end{rem}

Как следствие лемм \ref{lem:1}, \ref{lem:3}  и \ref{lem:6}
получаем следующий важный для дальнейшего результат.

\begin{lem}\label{lem:7}
Пусть число  \(\mu\in\mathbb C^-\) такое, что
\(\|G(\mu)\|<1/2\). Тогда  оператор  $A$ обладает инвариантным,  максимальнным
неотрицательным подпространством в том и только в том случае,
когда найдется нерастягивающий оператор $K: \mathcal H^+
\to\mathcal H^-$, такой, что
\begin{equation}\label{Ric2}
    F+K=(A_{22}-\mu)^{-1}(1-KG)^{-1}K(S-\mu).
\end{equation}
\end{lem}
Решение $K$ этого уравнения не зависит от $\mu\in\rho(A_{22})$.

\begin{proof}
 Первое утверждение леммы есть следствие лемм \ref{lem:1} и
 \ref{lem:6}. Далее, уравнение \eqref{Ric2} эквивалентно уравнению
 \eqref{Ric1}, где операторы в обеих частях равенств определены на
 $\mathcal D^+.$  Поскольку операторы $A_{21}=(A_{22} -\mu)F$ и
 $KA_{11}$ корректно опеределены на $\mathcal D^+$, то
 $(A_{22}-\mu)K$  и  $KA_{12}F$ обладают тем же свойством. Но
 тогда в \eqref{Ric1} можно раскрыть скобки и получить уравнение
 \eqref{eq:R2}. Тем самым, $K$  не зависит от $\mu$.
 \end{proof}

 Обозначим  \(\mathcal H_S=\mathcal
D(\overline{S})\subset\mathcal H^+\), где
\(\overline{S}\) --- замыкание оператора  \(S\). Снабженный нормой графика
\[
    \|x_+\|_{\mathcal H_S}=\sqrt{\|\overline{S}x_+\|^2+\|x_+\|^2},
\]
линеал $\mathcal H_S$ является гильбертовым пространством, плотно вложенным
в $\mathcal H^+$.
\begin{lem}\label{lem:8}
Существует ортонормированная система
\(\{\varphi_k\}_1^{\infty}\) в пространстве \(\mathcal H^+\),
которая с точностью до нормировки
является базисом Рисса в пространстве
\(\mathcal H_S\).
\end{lem}

\begin{proof}
Из определения нормы в пространстве $\mathcal H_S$ следует, что оператор
$C=1+(S^*S)^{1/2}$ изоморфно отображает $\mathcal H_S$ в $\mathcal H$. Если оператор
$C^{-1}$ компактный, а $\{\mu_k\}_1^{\infty}$ ---
последовательность его собственных значений, то в качестве искомой
системы можно взять полную ортонормированную в $\mathcal H$ систему
$\{\psi_k\}_1^{\infty}$ собственных векторов этого оператора.
Очевидно,  система $\{\mu_k^{-1}\psi_k\}_1^{\infty}$ будет полной
и ортонормированной в  пространстве $\mathcal H_S$.

Если $C^{-1}$ некомпактен, то построим  оператор $C_1\geqslant 1$,
такой, что $1\geqslant C_1-C\geqslant 0$, а спектр $C_1$ состоит
только из собственных значений. Если $E_{\lambda}$ ---
спектральная функция оператора $C$, то можно взять $C_1=
1+\sum_{k=1}^{\infty}k(E_{k+0} -E_{(k-1)})$. Выберем
ортонормированный базис $\{\psi_{k,j}\}_{j=1}^{\infty}$ в каждом
подпространстве $\mathcal H_k$, совпадающим с образом проектора
$E_{k+0} -E_{(k-1)}$. Тогда объединенная система
$\{\psi_{k,j}\}_{k,j=1}^{\infty}$, состоящая из собственых функций
оператора $C_1$, является полной и ортонормированной в $\mathcal
H$, а система $\{k\psi_{k,j}\}_{j=1}^{\infty}$, очевидно, является
эквивалентной ортонормированному базису в $\mathcal H_S$, т.е.
базисом Рисса.
\end{proof}

\begin{lem}\label{lem:9}
Пусть $\mathcal H_1$ и $\mathcal H_2$ гильбертовы пространства и последовательность
операторов $K_n: \mathcal H_1
\to \mathcal H_2$ слабо сходится к оператору $K: \mathcal H_1 \to \mathcal H_2$
(будем писать \(K_{n}\rightharpoonup K\)). Если $G: \mathcal H_2\to
\mathcal H_1$ --- компактный оператор, то последовательность $K_n\, G\, K_n$
сходится к $ K\, G\, K$ в равномерной операторной топологии (пишем
$K_n\, G\, K_n \Rightarrow K\, G\, K$).
\end{lem}
\begin{proof} Компактный оператор $G$ можно с произвольной
точностью приблизить конечномерным оператором. Поэтому достоточно
доказать утверждение леммы для одномерного оператора $G= (\cdot ,
\phi)\psi$, где $\phi\in H_1, \ \psi\in H_2$. Но для такого случая
утверждение сразу следует из определения слабой сходимости.
\end{proof}

\begin{proof}[Доказательство Теоремы~\ref{tm:1}.]
Обозначим через \(P_n\) ортогональные проекторы на $n$-мерные
подпространства $\Lin\{\varphi_k\}_1^n$, предполагая, что система
$\{\varphi_k\}_1^n$
обладает свойствами, указанными в Лемме~\ref{lem:8}. Тогда \(P_n\to 1\)
(сильно) в пространстве  \(\mathcal H^+\), а также
\(P_n\to 1\) в \(\mathcal H_S\). Рассмотрим  операторы
\[
    A_n=\begin{pmatrix}P_nA_{11}P_n&P_nA_{12}\\ A_{21}P_n&A_{22}
    \end{pmatrix},
\]
действующие в пространствах $\mathcal H_n^+\oplus \mathcal H^-$, где
    $ \mathcal H_n^+=P_n(\mathcal H^+)$.
Тогда  \(A_n\) являются \(m\)-диссипативными операторами в
пространстве Понтрягина \(\Pi_{\varkappa}\). Согласно теореме
Крейна-Лангера-Азизова операторы $A_n$ обладают $n$-мерными
максимальными неотрицательными подпространствами. В силу леммы
\ref{lem:7} найдутся нерастягивающие операторы $K_n: \mathcal H_n^+\to
\mathcal H^-$, такие, что
\begin{equation}\label{E1}
    F_n+K_n=(A_{22}-\mu)^{-1}(1-K_nG)^{-1}K_n(S_n-\mu).
\end{equation}
Известно, что единичный шар в сепарабельном гильбертовом
пространстве является слабо компактным. Так как последовательность
операторов $K_n$ такова, что их нормы не превосходят единицы, то
можно выбрать слабо сходящуюся подпоследовательность
\(K_{n_j}\rightharpoonup K\). Далее для краткости опускаем индекс $j$.
Конечно, норма предельного оператора $K$  также не превосходит
единицы. Учитывая, что $F$  ограниченный, а $G$ компактный
операторы и пользуясь леммой \ref{lem:9}, получаем
\[
    F_n=FP_n\to F, \quad K_nG\Rightarrow KG, \quad
    (1-K_nG)^{-1}\Rightarrow (1-KG)^{-1}.
\]
Имеем также
\begin{gather*}
    K_nS_n=K_nSP_n,\\
    \overline{S}P_nx\to\overline{S}x\qquad\forall x\in
    \mathcal D(\overline{S}).
\end{gather*}
Поэтому \(K_nSP_nx\rightharpoonup KSx\) и мы можем перейти к
слабому пределу в уравнении~\eqref{E1}. Тогда для предельного оператора \(K\)
получаем уравнение~\eqref{Ric2}.
В силу Леммы~\ref{lem:7} оператор $A$ имеет максимальное
неотрицательное инвариантное подпространство.

Докажем вторую часть теоремы о существовании инвариантного подпространства
$\mathcal L^+$, такого, что спектр сужения
\(A_+=\overline{A}|_{\mathcal L^+}\) лежит в замкнутой левой
полуплоскости.

Имеем
\begin{equation}\label{A}
    (\overline{A}-\mu)\begin{pmatrix}x_+\\ Kx_+\end{pmatrix}=
    \begin{pmatrix}(\overline{S}-\mu+GL)x_+\\ Lx_+\end{pmatrix},
\end{equation}
где \(L:=(A_{22}-\mu)(F+K)\), \(\mathcal D(L)=
\mathcal D(\overline S)\).

Обозначим через $Q$ ортопроектор из подпространства $\mathcal L$ на
$\mathcal H^+$, определенный формулой
\[
    Q\begin{pmatrix}
    x_+\\ Kx_+\end{pmatrix}=x_+ \qquad x= \begin{pmatrix}
    x_+\\ Kx_+\end{pmatrix} \in \mathcal L.
\]
Очевидно, \(Q\) ограниченно обратим, причем \(\|Q^{-1}\|\leqslant 2\).

Из равенств  ~\eqref{Ric1} и ~\eqref{A}  имеем

\begin{equation}\label{B}
    (\overline{A} - \alpha)|_{\mathcal L^+}=Q^{-1}(\overline S - \alpha +GL)Q=
    Q^{-1}[\overline S - \alpha +G(1-KG)^{-1}
    K(\overline{S}-\mu)]Q.
\end{equation}
Здесь мы используем Лемму~\ref{lem:3}  и предполагаем, что число
$\mu\in\mathbb C^-$ выбрано так, что $\|G(1-KG)^{-1}K\| < 1/2$.
Тогда при больших положительных $\alpha$ правая часть является
обратимым оператором в том и только том случае, когда оператор
$\overline S -\alpha$ обратим, т.е. когда $\overline S$ является
максимально диссипативным. Согласно Лемме~\ref{lem:3} это
эквивалентно максимальной диссипативности оператора $A$  в
пространстве $(\,
\mathcal H, [\, \cdot ,\cdot\, ]\, )$.
Следовательно, открытая левая полуплоскость может принадлежать
резольвентному множеству сужения $\overline{A}|_{\mathcal L^+}$
тогда и только тогда, когда оператор $A$ максимально диссипативный
в пространстве $(\, \mathcal H, [\, \cdot ,\cdot\, ]\, )$.

Пусть $A$ обладает указанным свойством. Тогда оператор $\overline
S -\alpha$ обратим при $\alpha\in\mathbb C^+$ и равенство~\eqref{B}
можно переписать в виде
\begin{gather}\label{E2}
    (\overline{A}-\alpha)|_{\mathcal L^+}=Q^{-1}[1+M(\alpha)]
    (\overline{S}(\mu)-\alpha)Q,\\ \intertext{где}\notag
    M(\alpha)=G(1-KG)^{-1}K[1+ (\alpha
    -\mu)(\overline{S}-\alpha)^{-1}]
\end{gather}
является голоморфной оператор-функцией, значения которой есть
компактные операторы.  При $\alpha\in\mathbb C^+$ имеем
$(\Re\alpha)\|(\overline S -\alpha)^{-1}\| \leqslant 1$. Полагая
$\alpha = -\ov\mu$ и пользуясь Леммой~\ref{lem:3}, получаем
$\|M(-ov\mu)\| \to 0$, если $\mu\to infty$  в секторе $|\arg (\pi-
\mu)|\leqslant \gamma <\pi/2$. В частности, число $\mu$ можно выбрать
так, чтобы оператор $1+M(\al)$ был обратим при $\al=-\mu$. Из
теоремы о голоморфной оператор-функции (см., например,
\cite[Гл. 1]{20} тогда следует, что спектр оператор-функции
$1+M(\alpha)$ в открытой правой полуплоскости дискретный (так как
значения $M(\alpha)$ "--- компактные операторы). Следовательно,
спектр сужения $A_+=\overline{A}|_{\mathcal L^+}$  в правой
полуплоскости может состоять только из изолированных собственных
значений конечной алгебрамческой кратности.

Заметим, что $\Im[Ax_0,x_0] =(\Im \alpha_0)[x_0,x_0]$, если
$Ax_0=\alpha_0 x_0$. Следовательно, если $A$  строго диссипативный
в пространстве $(\, \mathcal H, [\, \cdot ,\cdot\, ]\, )$ ( а
тогда его сужение $\overline{A}|_{\mathcal L^+}$ обладает тем же
свойством), то оператор $\overline{A}|_{\mathcal L^+}$ не имеет
собственных значений в открытой правой полуплоскости, а потому вся
эта полуплоскость принадлежит резольвентному множеству. Тем самым
теорема полностью доказана для строго диссипативных операторов
$A$.

В общем случае рассмотрим семейство операторов
\[
    A_{\varepsilon}=A-\varepsilon P_+,\qquad\varepsilon>0.
\]
Уже доказано, что для этих операторов имеются макксимальные,
неотрицательные, $A_{\varepsilon}$-инвариантные подпространства с
угловыми операторами $K_{\varepsilon}$. Сужения операторов
$A_{\varepsilon}$ на эти инвариантные подпространства являются
строго диссипативными операторами, так как
\[
    \Re\left[(A-\varepsilon P_+)\begin{pmatrix}x_+\\ Kx_+\end{pmatrix},
    \begin{pmatrix}x_+\\ Kx_+\end{pmatrix}\right]\leqslant
    -\varepsilon(x_+,x_+).
\]
Поэтому открытая правая полуплоскость лежит в резольвентных
множествах этих сужений.

Запишем модифицированное  уравнение Риккати для оператора
\(A_{\varepsilon}\):
\[
    F+K_{\varepsilon}=(A_{22}-\mu)^{-1}(1-K_{\varepsilon}G)^{-1}
    K_{\varepsilon}(S+i\varepsilon-\mu).
\]
Выберем подпоследовательность \(\varepsilon_n\to 0\)  так, чтобы
\(K_{\varepsilon_n}=:K_n\rightharpoonup K\). Сужения операторов \(A_{\varepsilon}\)
на соответсвующие инвариантные подпространства \(\mathcal L_{\varepsilon}^+\)
имеют представления
\begin{gather*}
    A_{\varepsilon}^+=A_{\varepsilon}|_{\mathcal L_{\varepsilon}^+}=
    Q^{-1}[1+M_{\varepsilon}(\alpha)](S+i\varepsilon-
    \alpha)Q,\\ \intertext{где}
    M_{\varepsilon}(\alpha)=G(1-K_{\varepsilon}G)^{-1}K_{\varepsilon}
    (S+i\varepsilon-\mu)(S+i\varepsilon-\alpha)^{-1}\Rightarrow
    M(\alpha).
\end{gather*}

Здесь мы вновь учитываем, что $(1-K_nG)^{-1}K_n \Rightarrow
(1-KG)^{-1}K$ при $\varepsilon =\varepsilon_n \to 0$, так как $G$
"--- компактный оператор. Оператор-функции
\(1+M_{\varepsilon}(\alpha)\) голоморфны в полуплоскости \(\mathbb
C_-\), и ограниченно обратимы в этой полуплоскости при любом
$\varepsilon>0$.  При $\varepsilon \to 0$ они  сходятся в
равномерной операторной топологии к оператор-функции
$1+M(\alpha)$, которая может иметь только изолированные
собственные значения в \(\mathbb C^-\). Если некая точка
$\alpha_0$ является ее собственным значением, то в малой
окрестности этой точки найдется собственное значение
\(1+M_{\varepsilon}(\alpha)\) при достаточно малом
$\varepsilon>0$, что не может быть. Этим завершается
доказательство теоремы для общего случая.

 \end{proof}

\section{Полугрупповые свойства сужений}

 Мы предполагаем, что читателю известны определения голоморфной и
$C_0$--полугрупп и их генераторов (см. например,
\cite[Ch. I, II]{22}, \cite[Ch. IX]{1}) в гильбертовом (или банаховом)
пространстве $\mathcal H$. Далее мы будем использовать только
свойства резольвент генераторов таких полугрупп. Напомним
следующие предложения. Первые два хорошо известны (см.
цитированные монографии), они дают критерии голоморфной и
$C_0$-полугруппы в терминах резольвенты генератора.

\begin{tmOG}[Feller-Mijadera-Phillips]
Оператор $T$ является генератором $C_0$--полугруппы тогда и только
тогда, когда найдётся число $\omega\in\mathbb{R}$, такое, что
полуплоскость $\Re\la>\omega$ лежит в резольвентном множестве
$\rho(T)$ и
\begin{equation}\label{3.1}
    \|(T-\la)^{-n}\|\le C(\la-\omega)^{-n}
\end{equation}
для всех целых $n\ge 1$ и $\la>\omega$, где постоянная $C$ не
зависит от $\la$ и $n$.
\end{tmOG}

В этом случае для полугруппы $U(t)=\exp(Tt)$ выполняется оценка
$\|U(t)\|\le C e^{\omega t}$. Нижняя грань чисел $\omega$, для
которых выполнена эта оценка или оценка \eqref{3.1}, называется
\textit{экспоненциальным типом} $C_0$--полугруппы. В случае
$\omega<0$ полугруппу называют \textit{экспоненциально
устойчивой}.

\begin{tmOGG}
Оператор $T$ является генератором голоморфной полугруппы тогда и
только тогда, когда правая полуплоскость лежит в резольвентном
множестве $\rho(T)$ и выполняется оценка
\begin{equation}
\label{3.2} \|(T-\la)^{-1}\|\le C|\la|^{-1}
\end{equation}
для всех $\la\colon\Re\la>0$, где постоянная $C$ не зависит от
$\la$.
\end{tmOGG}

Конечно, такая оценка выполняется, если $T$ является максимальным
секториальным оператором, т.е. его числовой образ ($=$ множество
значений квадратичной формы $(Tx,x)$, когда $x$ пробегает $D(T)$ и
$\|x\|=1$) лежит в секторе $\Lambda_\alpha=\{\la\colon
|\arg(\la-\pi)|<\alpha\}$ при некотором $\alpha<\pi/2$. В этом
случае неравенство \eqref{3.2} выполняется при
$C=(\cos\alpha)^{-1}$.

Далее удобно будет ввести следующее определение: полугруппу
$U(t)=\exp(Tt)$ назовём \textit{квази-голоморфной}, если для
резольвенты её генератора выполнена оценка \eqref{3.2} для всех
$\la\colon\Re\la>\ve$, где $\ve>0$ произвольно, а постоянная $C$
зависит только $\ve$. Очевидно, это условие эквивалентно тому, что
$T-\ve$ является генератором голоморфной полугруппы при любом
$\ve>0$. Известно, что квази-голоморфная полугруппа является
$C_0$-полугруппой экспоненциального типа 0. Обратное, конечно,
неверно.

Важную роль в дальнейшем играют следующие две теоремы.

\begin{tmGo}[\cite{21}]
Оператор $T$ является генератором $C_0$-полугруппы в пространстве
$\mathcal H$, допускающей оценку $\|U(t)\|\le C e^{\omega t}$,
тогда и только тогда, когда $\rho(T)$ содержит полуплоскость
$\Re\la>\omega$ и для любого $x\in \mathcal H$
\begin{equation}
\label{3.3} \sup\limits_{\delta>\omega} (\delta - \omega)
\int\limits_{\delta-i\infty}^{\delta+i\infty}\left (
\|(T-\la)^{-1}x\|^2+\|(T^\ast-\la)^{-1}x\|^2\right ) |d\la|<\infty
\end{equation}
\end{tmGo}

\begin{tmGe}[\mbox{\cite[Теорема V.1.11]{22}}]
Экспоненциальный тип $C_0$-полугруппы, генерируемой оператором $T$, совпадает с
нижней гранью чисел $\beta$, для которых полуплоскость
$\Re\la>\beta$ лежит в $\rho(T)$ и
\begin{equation}
\label{3.4} \sup\limits_{\Re\la>\beta}\|(T-\la)^{-1}\|<\infty.
\end{equation}
\end{tmGe}

\begin{lem}\label{lem:3.5}
Оператор $T$ генерирует  $C_0$--полугруппу экспоненциального типа
$0$ (экспоненциально устойчивую), если спектр $T$ лежит  в
замкнутой левой полуплоскости, оценка \eqref{3.3} выполняется при
достаточно больших $\omega>0$, а оценка \eqref{3.4} при всех
$\beta>0$ (при $\beta=0$).
\end{lem}

\begin{proof}
Утверждение следует из теорем Гомилко и Герхарта. Нужно только
заметить, что если оценка \eqref{3.4} выполняется при $\beta=0$,
то найдётся $\ve>0$ такое, что  $\rho(T)$ содержит полуплоскость
$\Re\la>-\ve$ и \eqref{3.4} выполняется при $\beta=-\ve$.
\end{proof}

\begin{lem}\label{lem:3.6}
Пусть оператор $-A_{22}$ является генератором квази-голоморфной
полугруппы и $G(\mu)= A_{12}(A_{22}-\mu)^{-1}$ "--- компактный
оператор при некотором $\mu=\mu_0\in\rho(A_{22})$. Тогда оператор
функция $G(\mu)$ корректно определена в левой  полуплоскости и
$$
\|G(\mu)\|\to 0 \text{ при } \mu\to\infty
$$
равномерно в полуплоскости $\Re\mu<0$. В частности, это свойство
выполняется, если $-A_{22}$ является секториальным оператором.
\end{lem}

\begin{proof}
Воспользуемся равенством
\[
    G(\mu)=G(\mu_0)(A_{22}-\mu_0)(A_{22}-\mu)^{-1}.
\]
Очевидно, в условиях леммы  норма оператор-функции
\((A_{22}-\mu_0)(A_{22}-\mu)^{-1}\) равномерно ограничена в
левой полуплоскости вне малой окрестности нуля.
 Компактный оператор
\(G(\mu_0):\mathcal H^-\to\mathcal H^+\) можно приблизить с
произвольной точностью в операторной норме конечномерным
оператором. Поэтому достаточно показать, что
\(\|Q(A_{22}-\mu_0)(A_{22}-\mu)^{-1}\|
\to 0\) при \(\mu\to\infty\) в левой полуплоскости
для любого одномерного оператора
\(Q=(\cdot,\varphi)\psi\), где \(\varphi\in\mathcal H^-, \psi \in\mathcal H^+\).
Заметим, что \(Q\) можно с произвольной точностью приблизить в
операторной норме оператором
\(Q_0=(\cdot,\varphi_0)\psi\), где \(\varphi_0\in\mathcal
D(A_{22}^*)\) (в условиях леммы сопряженный оператор $A^*_{22}$
плотно определён). Но оператор
\(Q_0(A_{22}-\mu_0)\) ограничен, а \(\|(A_{22}-\mu)^{-1}\|\leqslant
C|\mu|^{-1}\)  в левой полуплоскости вне малой окрестности нуля,
что следует из  теоремы о генераторе голоморфной полугруппы.
 Отсюда получаем утверждение леммы.
\end{proof}

\begin{lem}\label{lem:3.7}
Пусть $T$ "--- $m$-диссипативный, а $V$ "--- компактный операторы в
$\mathcal H$. Тогда при любом $\ve>0$ найдется достаточно большое
число $r=r(\ve)$, такое, что  для всех $\lambda$ из полуплоскости
$\Re \lambda\geqslant \ve$ и $|\lambda>r$ оператор $T+V-\lambda$
обратим и
\begin{equation}\label{T+V}
\|(T+V-\lambda)^{-1}\| \leqslant (2\ve)^{-1}.
\end{equation}
\end{lem}

\begin{proof}
Действуем также, как в Лемме \ref{lem:3.6}. Пусть  сначала
$V=(\cdot ,\varphi )\psi$ "--- одномерный оператор в $\mathcal H$.
Представим его в виде
\[
    V=V_0+V_1,\qquad V_0=(\cdot,\varphi_0)\psi,\qquad
    \|V_1\|<\delta,\qquad \phi_0 \in \mathcal D(T^*).
\]
Тогда
\[
    (T-\lambda)^{-1}V=(T-\lambda)^{-1}V_1-\lambda^{-1}(V_0-(T-\lambda)^{-1}
    TV_0).
\]
Фиксируем \(\varepsilon>0\) и возьмём \(\delta<\varepsilon/4\).
Оператор \(SV_0\) ограничен, поэтому при $\Re \lambda\geqslant\ve$
имеем оценку
\[
    \|(T-\lambda^{-1}V\|\leqslant\dfrac{\delta}{\varepsilon}+
    \dfrac{\mathrm{const}}{|\lambda|\varepsilon}<\dfrac{1}{2},
\qquad\text{если}\quad |\lambda|>\frac{4\, const}\ve.
\]
Но тогда для таких $\lambda$ справедлива оценка \eqref{T+V}.
Доказательство сохраняется, если $V$ "--- конечномерный оператор.
Так как компактный оператор можно приблизить по норме  с
произвольной точностью конечномерным, то оценка остается
справедливой и для компактного оператора $V$.
\end{proof}

Заметим, что для максимально диссипативного оператора $-A_{22}-1$
корректно определены дробные степени (см. например, \cite[Гл.II]
{22}).

\begin{lem}\label{lem:3.8}
Пусть оператор $-A_{22}$ является генератором квази-голоморфной
полугруппы, а $A_{12}$ и $A_{21}$ таковы, что операторы
\begin{equation} \label{Pair}
A_{12}(-A_{22}-1)^\alpha \quad \text{и}\quad
(-A_{22}-1)^{1-\alpha}A_{21}
\end{equation}
ограничены при некотором $0\le\alpha\le 1$. Тогда оператор-функция
$$
R(\mu)\colon=A_{12}(A_{22}-\mu)^{-1}A_{21}
$$
равномерно ограничена в полуплоскости $\Re\mu<0$ вне малой
окрестности нуля.
\end{lem}

\begin{proof}
Из условия леммы следует, что оператор $R(-1)$ ограничен. Далее,
$$
R(\mu)=R(-1)-A_{12}(-A_{22}-1)^\alpha\left [
(\mu+1)(A_{22}-\mu)^{-1} \right ](-A_{22}-1)^{1-\alpha}A_{21}.
$$
Из теоремы о генераторе голоморфной полугруппы следует  оценка
$\|(\mu+1)(A_{22}-\mu)^{-1}\|\le C$ для всех $\mu$  в левой
полуплоскости вне окрестности нуля, из которой получаем
утверждение леммы.
\end{proof}

Отметим, что для оператор-матриц, возникающих в конкретных задачах
(см., например,~\cite{14,15}), ограниченность операторов
\eqref{Pair} бывает несложно проверить при $\al=1/2$.

\begin{lem}\label{lem:3.9}
Пусть $A_+$ --- сужение оператора $\overline{A}$ на  инвариантное,
максимальное неотрицательное подпространство $\mathcal{L}^+$.
Тогда
\begin{equation}
\label{A_+} A_+=Q^{-1}XQ,\quad
X=\overline{S}+G(1-KG)^{-1}K(\overline{S}-\mu),
\end{equation}
а сопряжённый оператор $A_+^*$ имеет представление
$$
A_+^*=Q^{-1}(1+K^*K)^{-1}X^*(1+K^*K)Q.
$$
Здесь $Q$ --- ортопроектор из $\mathcal{L}^+$ на $\mathcal{H}^+$.
Операторы $K$ и $X$ не зависят от $\mu$.
\end{lem}

\begin{proof}
Представление \eqref{A_+} вытекает из (9). Пусть $A_+^*=Q^{-1}YQ$.
Тогда
\begin{multline*}
\left ( A_+
\begin{pmatrix}
x\\Kx
\end{pmatrix},
\begin{pmatrix}
y\\Ky
\end{pmatrix}
\right ) = \left (
\begin{pmatrix}
Xx\\KXx
\end{pmatrix},
\begin{pmatrix}
y\\Ky
\end{pmatrix}
\right )=\\
=\left ( (1+K^*K)Xx,y \right )= \left (
\begin{pmatrix}
x\\Kx
\end{pmatrix},
A_+^*\begin{pmatrix} y\\Ky
\end{pmatrix}
\right )= (x,(1+K^*K)Yy),
\end{multline*}
что влечёт $X^*(1+K^*K)=(1+K^*K)Y$. В силу леммы \ref{lem:7}
оператор $K$ (а потому и $X$) не зависит от $\mu$.
\end{proof}

Теперь докажем  теоремы о полугрупповых свойствах оператора $A_+$.

\begin{tm}
Пусть  выполнены условия теоремы 2.1 и $A^+$ "--- сужение
$m$-диссипативного оператора $\overline A$ в пространстве
$(\mathcal{H},[\cdot,\cdot])$ на инвариантное, максимальное
неотрицательное подпространство $\mathcal L^+$, такое, что $\sigma
(A^+)$ лежит в замкнутой левой полуплоскости. Пусть также
выполнено одно из условий
\begin{itemize}
\item[(i)] оператор \(A_{12}\) компактен;
\item[(ii)] при любом $\ve>0$ оператор-функция
$R(\mu)=A_{12}(A_{22}-\mu)^{-1}A_{21}$ равномерно ограничена в
полуплоскости $\mathbb{C}_\ve^-=\{\mu\colon\Re\mu\le\ve\}$, а
$\|G(\mu)\|\to 0$ при $\mu\to\infty$ равномерно в
$\mathbb{C}_\ve^-$;
\item[(iii)] оператор $-A_{22}$ генерирует квази-голоморфную
полугруппу и оба оператора в \eqref{Pair} ограничены при некотором
$0\le\alpha\le 1$.
\end{itemize}
Тогда $A_+$ генерирует $C_0$-полугруппу экспоненциального типа
$0$. Если $\overline{A}$ --- равномерно  $m$-диссипативный в
пространстве $(\mathcal{H},[\cdot,\cdot])$, т.е. $\Re[Ax,x]\le
-\ve(x,x)$, $\ve>0$, для всех  $x\in D(A)$, то $A_+$ генерирует
экспоненциально устойчивую $C_0$-полугруппу.
\end{tm}

\begin{proof}
Пусть выполнено условие~(i). Из формулы~\eqref{A} имеем
\[
    A_+=Q^{-1}(\overline{S}+GL)Q=Q^{-1}(\overline{S}+A_{12}(F+K))Q.
\]
Оператор \(\overline{S}=\overline{S}(\mu)\) является
$m$-диссипативным, а потому генерирует полугруппу сжатий в
\(\mathcal H_+\), тем более, \(C_0\)-полугруппу. Известно, что
ограниченное возмущение такого оператора (см.~\cite[Гл.~5]{14})
является генератором
\(C_0\)-полугруппы.

Докажем, что экспоненциальный тип этой полугруппы равен нулю. Мы
знаем, что спектр $A^+$ лежит в замкнутой левой полуплоскости. Из
Леммы \ref{lem:3.7} тогда следует, что резольвента
$(\overline{S}+A_{12}(F+K))^{-1}$ является ограниченной в
полуплоскости $\Re \alpha >\ve$  при любом $\ve>0$. Остается
воспользоваться теоремой Герхарта.

Докажем теперь утверждение теоремы при выполнении условия (ii).
Фиксируем число $\ve>0$. Для оператора $A_+$ воспользуемся
представлением \eqref{A_+},  из которого получим

\begin{multline}\label{X}
(A_+-\alpha)=Q^{-1}(X-\alpha)Q= Q^{-1}((\ov
S(\mu)-\alpha)[1+M_1(\mu, \alpha)]Q,\\
 M_1(\mu, \alpha)= [1+(\alpha -\mu)
(\overline{S}-\alpha)^{-1}]G(1-KG)^{-1}K,
\end{multline}
где $G$ и $S$ зависят от $\mu$. Положим $\mu =-\ov\alpha$  и
выберем число $r_0$ выбрано столь большим, что
$\|G(\overline{-\alpha})\|<1/5$ при всех $|\al|>r_0, \Re\al >0$.
Так как оператор $\overline{S}$ является $m$-диссипативным, то
$|\alpha-\ov\alpha|\cdot|(\overline{S}-\alpha)^{-1}\|\le 2$,
поэтому
$$
\|M(-\ov\alpha, \alpha)\|< 2\cdot(1/5)\cdot(5/4) = 1/2.
$$
Учитывая, что $\|Q^{-1}\|\le 2$  и $\|Q\|\le 1$, для всех
$\alpha\colon\Re\alpha\ge\ve$, $|\alpha|>r_0$ получаем
$$
\|(A_+-\al)^{-1}\|\le 4\|(\ov{S}(-\ov\al)-\al)^{-1}\|\le 4\ve\p .
$$
Поскольку открытая правая полуплоскость принадлежит $\rho(A_+)$,
то $(A_+-\al)\p$ ограничена в полуплоскости
$\al\colon\Re\al\ge\ve$ при любом $\ve>0$.

Учитывая Леммы \ref{lem:3.5} и \ref{lem:3.9}, остаётся показать
что при больших $\omega$   и любом $x_+\in \mathcal H^+$
справедливы оценки
\begin{equation}
\label{GoX}
\sup\limits_{\delta>\omega}(\delta -\omega)\int\limits_{\delta-i\infty}^{\delta+i\infty}
\|(X -\al)\p x\|^2 \, |d\al|, \quad
\sup\limits_{\delta>\omega}(\delta -\omega)\int\limits_{\delta-i\infty}^{\delta+i\infty}
\|(X^* -\al)\p x\|^2 \, |d\al|.
\end{equation}
где
$$
X-\alpha = (\ov S -\alpha)(1+ M_1(-\ov\alpha, \alpha)).
$$

 Уже показано, что при $\Re\al > r_0$ выполнена оценка
$\|M(-\ov\al, \al)\|< 1/2$. Поэтому первую оценку достаточно
доказать, если в ней вместо $X$  участвует $\ov S =\ov
S(-\ov\al)$.

Спараведливо равенство
$$
A_{11}-\al=\ov S (-\ov\al)-\al-R(-\ov\al),
$$
и по условию $\|R(-\ov{\al})\|\le C$
 при $\al\colon\Re\al\ge\ve$. Тогда из
 оценки $\|(\ov S(-\ov\al)-\al)\p\|\le 1/\Re\al$ следует, что
$A_{11}-\al$ обратим при $\Re\al\ge 2C$ (а потому обратим в правой
полуплоскости, т.к. $A_{11}$ диссипативный). Кроме того, при
$\Re\al\ge 2C$ справедливы оценки
$$
2/3\le\|(1+(A_{11}-\al)\p R(-\ov\al))\p\|<2.
$$
Поэтому при $\omega>\max\{r_0,2C\}$  первая оценка в
\eqref{GoX} выполняется тогда и только тогда, когда она выполняется
после замены $X$ на  $A_{11}$. Но $A_{11}$ является
$m$-диссипативным и генерирует полугруппу сжатий. Поэтому такая
оценка для $A_{11}$ заведомо выполнена.

Для доказательства второй оценки надо представить оператор $X$ в
виде
$$
X-\al = (1+M(\mu,\al))(\ov S -\al), \qquad M(\mu,\al)=
G(1-GK)^{-1}\left(1+(\alpha-\mu)(\ov S -\al)^{-1}\right).
$$
 и повторить предыдущие рассуждения.
 Тем самым, доказана достаточность условия (ii).

Если выполнено условие (iii), то в силу Лемм \ref{lem:3.6} и \ref{lem:3.7}
выполнено условие (ii).

Пусть теперь оператор $\ov A$ равномерно $m$-диссипативный
пространстве $(\,\mathcal H, [\cdot , \cdot]\,)$. Из Леммы~\ref{lem:5} тогда
следует, что при некотором $\ve>0$
$$
\Re (S(\mu)x_+,x_+)\le -\ve(x_+,x_+),\quad x_+\in D^+,
$$
для всех $\mu\in\mathbb{C}^-$. Тогда число $\mu\in\mathbb{C}^-$
можно выбрать так, чтобы $\|M_1(\mu, \alpha)\|<1/2$ для всех $\al$
из правой полуплоскости. В этом случае, учитывая равномерную
диссипативность оператора $S(\mu)$, получаем
$$
\|(A_+-\al)\p\|\le 4\|(\ov S(\mu)-\al)\p\|\le 4\ve\p
$$
при всех $\Re\al>0$. Согласно Лемме~\ref{lem:3.9} оператор $A_+$
генерирует экспоненциально устойчивую $C_0$-полугруппу. Теорема
полностью доказана.
\end{proof}

\begin{tm}
Пусть  выполнены условия теоремы 2.1 и $A^+$ "--- сужение
$m$-диссипативного оператора $\overline A$ в пространстве
$(\mathcal{H},[\cdot,\cdot])$ на инвариантное, максимальное
неотрицательное подпространство $\mathcal L^+$, такое, что $\sigma
(A^+)$ лежит в замкнутой левой полуплоскости. Пусть также
выполнено одно из условий

\begin{itemize}
\item[(i)] оператор $\ov S=\ov S(\mu)$ генерирует квази-голоморфную
полугруппу при некотором $\mu\in\mathbb C^-$;
\item[(ii)] оператор $A_{11}$ генерирует
квази-голоморфные полугруппу, причём оператор
$A_{21}(A_{11}-\al_0)\p$ ограничен при некотором $\al_0\in\mathbb
C^+$.
\end{itemize}
Тогда $A_+$ генерирует квази-голоморфную полугруппу.
\end{tm}
\begin{proof}
Пусть \(\overline{S}(\mu_0)\) генерирует квази-голоморфную полугруппу.
Из формулы
\[
    R(\mu)=A_{12}(A_{22}-\mu)^{-1}A_{21}=R(\mu_0)+(\mu-\mu_0)
    G(\mu)F(\mu_0)
\]
следует, что при любом \(\mu\in\mathbb C^-\) оператор \(S(\mu)\) является
ограниченным возмущением \(S(\mu_0)\). Кроме того, \(\overline{S}(\mu)\)
является \(m\)-диссипативным, а потому он также генерирует квази-голоморфную
полугруппу.

Воспользуемся равенством
$$
\ov S(\mu)-\al = [1-(\mu-\mu_0)G(\mu)F(\mu_0)(\ov S_0-\al)^{-1}](\ov S_0-\al),
$$
где $\ov S_0=\ov S(\mu_0).$ Найдется достаточно большое число
$r=r(\mu,\ve)$, такое, что при $|\al|>r,$  $\Re\al \geqslant\ve$
выполняется оценка
$$
\|(\mu-\mu_0)G(\mu)F(\mu_0)(\ov S_0-\al)^{-1}\| <1/2.
$$
Слндовательно,
\begin{equation}\label{S}
\|(\ov S(\mu)-\al)^{-1}\|<2\|(\ov S_0-\al)^{-1}\|<C(\ve)|\al|^{-1}
\quad\text{при}\  |\al|>r(\mu,\ve),\ \, \Re\al\geqslant\ve.
\end{equation}
Можно считать, что $C=C(\ve)\geqslant 1$. Фиксируем $\ve>0$  и
выберем $\mu<0$ так, чтобы
$$
\|G(\mu)(1-KG(\mu))^{-1}\| < (4C)^{-1}.
$$
Тогда для функции $M_1(\mu,\al)$, определенной в \eqref{X},
получаем оценку
$$
\|M_1(\mu,\al)\|<(4C)^{-1}(1+C(1+|\mu|\, |\al|^{-1}) <3/4,
$$
если $|\al|>\max\{r,|\mu|\}$. Так как $\ov S(\mu)$ генерирует
квази-голоморфную полугруппу, то из представления \eqref{X} имеем
\begin{equation}\label{est}
    \|(X-\alpha)^{-1}\|\leqslant 4\|(\overline{S}(\mu)-\alpha)^{-1}\|
    \leqslant C|\alpha|^{-1}\qquad
    \text{при }\Re\alpha\geqslant\ve, |\al|>r(\mu,\ve),
\end{equation}
с постоянной $C$, не зависящей от $\al$.  Так как спектр \(X\)
лежит в левой полуплоскости, то эта оценка остаётся справедливой
во всей  полуплоскости
\(\Re\alpha\geqslant\varepsilon\).
 Следовательно, \(X\) генерирует квази-голоморфную полугруппу.

Покажем, что условие~(ii) влечет условие~(i). Воспользуемся
равенством
\[
    S(\mu)-\alpha=A_{11}-\alpha-R(\mu)=
    \left[1-G(\mu)A_{12}(A_{11}-\alpha)^{-1}\right]\cdot(A_{11}-\alpha).
\]
Поскольку $A_{11}$ генерирует квази-голоморфную полугруппу, то
при всех \(\alpha\) из полуплоскости
\(\Re\alpha\geqslant\varepsilon\) справедлива оценка
\[
    \|A_{12}(A_{11}-\al)^{-1}\|\leqslant
    \|A_{21}(A_{11}-\alpha_0)^{-1}\|\,\|(A_{11}-\alpha_0)
    (A_{11}-\alpha)^{-1}\|\leqslant C(\varepsilon).
\]
Выберем \(\mu<0\) так, чтобы \(\|G(\mu)\|<1/2C(\varepsilon)\).
Тогда
\[
    \|(S(\mu)-\alpha)^{-1}\|\leqslant
    2\|(A_{11}-\alpha)^{-1}\|<C\,|\alpha|^{-1},
    \qquad\text{если} \ \Re\al\geqslant\, \ve.
\]
Следовательно, при выбранном значении $\mu$ оператор $S(\mu)$
генерирует квази-голоморфную полугруппу. Теорема доказана.
\end{proof}

\end{document}